\documentclass{article}
\usepackage[ansinew]{inputenc}
\usepackage[english]{babel}
\usepackage{graphicx}
\usepackage{amsmath,amssymb,amsthm,mathrsfs,amsfonts}
\usepackage{cancel}

\newtheorem {thm}{Theorem}[section]
\newtheorem {defn}{Definition}[section]

\newtheorem {lemma}[thm]{Lemma}
\newtheorem {problem}[defn]{Problem}

\newtheorem {cor}[thm]{Corollary}

\newcommand{\I}{\boldsymbol{\mathcal{I}}}

\newcommand{\U}{\boldsymbol{\mathcal{U}}}

\newcommand{\GP}{\boldsymbol{\mathcal{GP}}}
\newcommand{\SP}{{\operatorname{sp}}}
\newcommand{\TR}{{\operatorname{tr}}}
\newcommand{\SPEC}{{\operatorname{spec}}}

\newcommand{\bbmatrix}[1]{\left[ \begin{array}{cccccccccccccccccc} #1 \end{array} \right]}

\title{Construction of a Sturm-Liouville vessel using Gelfand-Levitan theory. On
solution of the Korteweg-de Vries equation in the first quadrant.}
\author{Andrey Melnikov\\Drexel University, Philadelphia, USA}
\begin{document}
\maketitle

\abstract{Using Gelfand-Levitan theory on a half line, we construct a vessel for the class of potentials, whose spectral functions
satisfy a certain regularity assumption. When the singular part of the spectral measure is absent, we construct a
canonical model of the vessel. Finally, evolving the constructed vessel, we solve the Korteweg de Vries equation on the half
line, coinciding with the given potential for $t=0$. It is shown that the initial value for $x=0$ is prescribed by this construction,
but can be perturbed using an ``orthogonal'' to the problem measure.

The results, presented in this work 1. include formulas for the ingredients of
the Gelfand-Levitan equation, 2. are shown to be general in the sense that NLS, Canonical systems and many more equations can be solved
using theory of vessels, analogously to Zacharov-Shabath scheme, 3. present a generalized inverse scattering theory on a line for potentials with singularities using pre-vessels,
4. present the tau function and its role.
}

\tableofcontents

\section{Introduction}
The Korteweg-de Vries (KdV) is the following nonlinear evolutionary Partial Differential Equation (PDE) for a
function of two real variables $q(x,t)$:
\begin{equation} \label{eq:KdV}
q_t = - \dfrac{3}{2} q q_x + \dfrac{1}{4} q_{xxx},
\end{equation}
where $q_t, q_x$ denote the partial derivatives.
The mathematical theory behind the KdV equation is rich and interesting, and, in the broad sense, is a topic of active mathematical research. The equation is named after Diederik Korteweg and Gustav de Vries who studied it in \cite{bib:KdV}.
Usually, one considers the \textit{initial value problem}, which is defined as follows: find a solution $q(x,t)$ of \eqref{eq:KdV}, which additionally satisfies:
\begin{itemize}
	\item $q(x,0) = q(x) (x\in\mathbb R)$, for a given function $q(x)$ and is called KdV on line, or
	\item $q(x,0) = q(x) (0\leq x\leq\infty)$, $q(0,t) = h(t) (0\leq t\leq\infty)$, called KdV in the first quadrant.
\end{itemize}
A more complicated problem consists of specifying the initial value on a closed loop.
A standard technique to solve KdV involves a more elementary equation, called Sturm Liouville (SL) differential equation:
\begin{equation}\label{eq:SL}
-\frac{d^2}{dx^2} y(x) + q(x) y(x) = \lambda y(x),
\end{equation}
where $\lambda\in\mathbb C$ is called the \textit{spectral parameter} and the coefficient $q(x)$ is called the \textit{potential}. In order to solve
\eqref{eq:KdV} using \eqref{eq:SL} one transforms \cite{bib:GGKM} the potential $q(x)$ appearing in \eqref{eq:SL} to 
its ``scattering data". Then one evolves with $t$ the scattering data using some simple differential equations.
Finally, transforming back the evolved scattering data we obtain a new potential $q(x,t)$ of two variables, which solves
\eqref{eq:KdV} and satisfies $q(x,0)=q(x)$. In other words, in this manner we
solve the initial value problem for the equation \eqref{eq:KdV}. So, in order to solve \eqref{eq:KdV} one
have to find "scattering data" for a given potential.

On the half line the question of characterizing of scattering data (or more precisely spectral function) for
a given potential was completely solved for a continuously differentiable potential by Gelfand-Levitan theory
\cite{bib:GL}, but not always it can be used to solve the KdV equation \eqref{eq:KdV}.
In the case the spectral function is ``differentiable''
($var[\rho(\lambda]<\infty$) it is possible to use it and this fact is reflected in many known solutions corresponding to this case:
\begin{enumerate}
	\item Soliton solutions correspond to $d\rho(\lambda)$ to be a finite sum of point mass measures (discrete measure),
	\item Krichever solution \cite{bib:Krich77}, where $d\rho(\lambda)=f(\lambda)d\lambda$ for $\lambda\in\Gamma$ 
		for some algebraic curve $\Gamma$,
	\item Fadeyev inverse scattering theory \cite{bib:FaddeyevII} where the $d\rho(\lambda)$ is supported on 
		the positive real line and has a finite number of point-mass measures on the negative real line,
	\item Periodic potentials \cite{bib:MagWin} 
		correspond to discrete spectral measures with accumulation point at infinity,
\end{enumerate}
and many more. In all these examples, the potential $q(x)$ exists on the whole line $\mathbb R$. Gelfand-Levitan
theory \cite{bib:GL}, on the other hand, was developed primary to answer the question of reconstructing the potential $q(x)$
from its spectral function $\rho(\lambda)$ \textit{on the half line} using some works of Marchenko and H. Weyl.

An analogous scheme for solution of the KdV equation \eqref{eq:KdV} was developed by the author using a theory of vessels
\cite{bib:KdVVessels}.
In this theory the scattering data is encoded (and easily translated) to a matrix-valued $2\times 2$ 
function $S(\lambda)$, which is analytic at infinity, with value $I$ there and is symmetric (see Theorem \ref{thm:SSym}).
As an initial step so called ``regular'' vessels were used, which involved only bounded operators. 
It is known \cite{bib:bgr} that such functions possess realizations as follows \eqref{eq:S0realized}
\[
S(\lambda) = I - B^*_0 \mathbb X_0^{-1}(\lambda I - A)^{-1}B_0\sigma_1, \quad \sigma_1=\bbmatrix{0&1\\1&0}, \quad I= id,
\]
where for an auxiliary Hilbert space $\mathcal H$ the operators act as follows: $\mathbb X_0, A: \mathcal H\rightarrow\mathcal H$, \linebreak $\quad B_0:\mathbb C^2\rightarrow\mathcal H$
and satisfy the following conditions
\[
A \mathbb X_0 + \mathbb X_0 A^* + B_0\sigma_1 B_0 = 0, \quad \mathbb X_0^* = \mathbb X_0.
\]
Evolving the operators with respect to $x$, one uniquely constructs \cite{bib:GenVessel} a vessel $\mathfrak{V}$, 
whose transfer function $S(\lambda,x)$
is realized as follows \eqref{eq:Srealized}:
\[ S(\lambda,x) = I - B^*(x) \mathbb X(x)^{-1}(\lambda I - A)^{-1}B(x)\sigma_1
\]
and coincides with the given $S(\lambda)$ for $x=x_0$ (see Definition \ref{def:Vessel} for details).
The operator $\mathbb X_0$ is assumed to be invertible and $A$ is assumed to be a generator of a $C_0$-group.
Notice that the singularities of $S(\lambda)$ are determined by the spectrum of $A$ only. Moreover, the same singularities with respect to $\lambda$ appear in $S(\lambda,x)$ and are $x$-independent. The new function
$S(\lambda,x)$ serves as a B\" acklund transformation (see Theorem \ref{thm:Backlund}) 
between trivial SL equation (having the zero potential)
and a more complicated one, having potential $2\dfrac{d^2}{dx^2} \ln(\det(\mathbb X(x)))$.  Evolving next
the operators with respect to $t$ as in Definition \ref{def:KdVVessel} we will obtain a KdV vessel, for which
$q(x,t) = 2\dfrac{d^2}{dx^2} \ln(\det(\mathbb X(x,t)))$ satisfies \eqref{eq:KdV}. As we can see this scheme is
essentially the same as the Zacharov-Shabbath's one \cite{bib:ZakSha74}, but has no assumption on the (rapid) decrease of the potential at infinity.
Moreover, one can use vessels to construct
solutions for other completely integrable PDEs, for example NLS was done in \cite{bib:ENLS}.

In this work we present a generalization of regular vessels, involving an unbounded operator $A$.
The main problem in this approach is that many formulas, present in a regular vessel, have to be rewritten, either
including the domain of $A$, or using its resolvent $R(\lambda) = (\lambda I - A)^{-1}$. We have chosen the second approach, so
for example, the Lyapunov equation
\[A \mathbb X_0 + \mathbb X_0 A^* + B_0\sigma_1 B_0 = 0
\]
becomes in the new setting \eqref{eq:LyapunovAt0}
\[ \mathbb X_0 R^*(-\bar\lambda) + R(\lambda) \mathbb X_0 - R(\lambda) B_0 \sigma_1 B_0^* R^*(-\bar\lambda) = 0,\quad
\lambda\not\in\SPEC(A). \]
If all the operators are bounded these two equations are equivalent.
All the regularity assumptions and formulas are rewritten so that Gelfand-Levitan theory fits into the setting
of vessels in a particular case of the spectral measure. More precisely, we construct a vessel 
in the case when a new function $w(\lambda)$, constructed from the spectral measure $\rho(\lambda)$
\[ \omega(\lambda) = \left\{ \begin{array}{ll}
\rho(\lambda) - \dfrac{2}{\pi}\sqrt{\lambda}, & \lambda > 0, \\
\rho(\lambda), & \lambda<0,
\end{array}\right.
\]
is a measure and has the property that \eqref{eq:Deff'} holds:
\[ f(x,y) = \dfrac{\partial^2 F}{\partial x\partial y} = \int\limits_{-\infty}^{+\infty} \cos(\sqrt{\mu}\,x)\cos(\sqrt{\mu}\,y) d\omega(\mu).
\]
The function $f(x,y)$ is one of the ingridients of the famous Gelfand-Levitan equation
\eqref{eq:GelfandLevitan}. We show in this case how to construct the corresponding KdV vessel, and as a result
a solution of the KdV equation \eqref{eq:KdV} coinciding with the given continuous $q(x)$ (whose spectral function
creates $f(x,y)$ via \eqref{eq:Deff'}). The value of the vessel for $x=0, t>0$ turns to be prescribed and we
present a formula for it \eqref{eq:q0t}
\[ q(0,t) = 2 \int_{\mathbb R} \sqrt{\lambda}\sin((\sqrt{\lambda})^3t) d\omega(\lambda). \]
It is also possible to perturb this value, without changing the value for $q(0,x)$, which is explained
in Section \ref{sec:KdV}.

Notice that the solution of the KdV equation $q(x,t) = 2\dfrac{d^2}{dx^2} \ln(\det(\mathbb X(x,t)))$ has
the property that it does not exist at the points where $\mathbb X(x,t)$ is singular. At these points
the potential $q(x,t)$ has a singularity, of course. So, we propose to study potentials, which arise in such
a manner using a notion of \textit{pre-vessel} (see Definition \ref{def:PreVessel}) .

At the last Section \ref{sec:supres} we show some models of vessels with an unbounded $A$ and a realization of the Fadeyev theory.
The canonical systems and Non Linear Schr\" odinger equations are presented as particular cases of the theory of
vessels as a final example. 
\section{Background}
\subsection{Summary of the Gelfand Levitan theory}
A fundamental paper of Israel Gelfand and Boris Levitan \cite{bib:GL} serves as a background and a reference for many
research works on SL equation \eqref{eq:SL} on a half line, 
.
Let $q(x) (0\leq x \leq \infty)$ be a continuous on each finite interval function. Consider the equation \eqref{eq:SL} with boundary condition
\begin{equation} \label{eq:SLinit}
y(0)=1,\quad y'(0)=h
\end{equation}
and let $\phi(x,\lambda)$ be its solution. There exists a representation of $\phi(x,\lambda)$ in the following form
\begin{equation} \label{eq:phicos}
\phi(x,\lambda) = \cos(\sqrt{\lambda}\,x) + \int\limits_0^x K(x,t) \cos(\sqrt{\lambda}\,t) dt.
\end{equation}
And conversely, there is a representation
\begin{equation} \label{eq:cosphi}
\cos(\sqrt{\lambda}\,x) = \phi(x,\lambda) - \int\limits_0^x K_1(x,t) \phi(t,\lambda) dt.
\end{equation}
From the work of Marchenko \cite{bib:MarchenkoSL} using Riemannian transformations, it follows that for continuously differentiable
$q(x)$ and $h=0$
\[ \begin{array}{lll}
\dfrac{\partial^2}{\partial x^2} K(x,t) - q(x) K(x,t) = \dfrac{\partial^2}{\partial t^2} K(x,t), \\
\dfrac{\partial}{\partial t} K(x,t)|_{t=0}, \quad \dfrac{d}{dx} K(x,x) = \int_0^x \dfrac{q(x)}{2} dx.
\end{array} \]
The second kernel satisfies
\begin{equation} \label{eq:PDEK1}
 \begin{array}{lll}
\dfrac{\partial^2}{\partial x^2} K_1(x,t) = \dfrac{\partial^2}{\partial t^2} K_1(x,t) - q(t) K(x,t), \\
\dfrac{\partial}{\partial t} K(x,t)|_{t=0}, \quad \dfrac{d}{dx} K(x,x) = \int_0^x \dfrac{q(x)}{2} dx.
\end{array} \end{equation}
These PDEs with initial value problems are known to have unique solutions.

There exists at least one non decreasing function $\rho(\lambda) (-\infty\leq\lambda\leq+\infty)$, such that for each
$f(x)\in L_2(0,\infty)$ functions
\[ E_n(\lambda) = \int\limits_0^n f(x) \phi(x,\lambda) dx
\]
converge in norm with respect to $\rho(\lambda)$ to a function $E(\lambda)$, i.e.
\[ \lim\limits_{n\rightarrow\infty} \int\limits_{-\infty}^{+\infty} [E(\lambda)-E_n(\lambda)]d\rho(\lambda)=0.
\]
In this case the Parseval equality holds:
\[ \int\limits_{-\infty}^{+\infty} E^2(\lambda) d\rho(\lambda) = \int_0^{+\infty} f^2(x)dx.
\]
The function $\rho(\lambda)$ is called the \textit{spectral function} of \eqref{eq:SL} with initial conditions
\eqref{eq:SLinit}. It is proved in Gelfand-Levitan paper that $\rho(\lambda)$ satisfies the conditions of the next Theorem.
\begin{thm}[\cite{bib:GL}]\label{thm:GLMnThm} Let $\rho(\lambda)$ be a monotone function. Let
\[ \omega(\lambda) = \left\{ \begin{array}{ll}
\rho(\lambda) - \dfrac{2}{\pi}\sqrt{\lambda}, & \lambda > 0, \\
\rho(\lambda), & \lambda<0
\end{array}\right.
\]
and suppose that $w(\lambda)$ satisfies
\begin{enumerate}
	\item[$1^\circ$] for all $x>0$ there exists integral 
		\[ \int\limits_{-\infty}^0 e^{\sqrt{|\lambda|}x} d\omega(\lambda),
		\]
	\item[$2^\circ$] the function
		\[ a(x) = \int\limits_1^\infty \dfrac{\cos(\sqrt{\lambda}\,x)}{\lambda}d\omega(\lambda)
		\]
	possesses the fourth continuous derivative.
\end{enumerate}
If the function $\rho(\lambda)$ (or more precisely $\omega(\lambda)$) satisfies conditions $1^\circ, 2^\circ$, then there exists a continuous function $q(x)$ and $h$, such that
$\rho(\lambda)$ is the spectral function of \eqref{eq:SL} with initial condition \eqref{eq:SLinit}. Conversely, if $q(x)$ is continuously
differentiable, the corresponding to it spectral function $\rho(\lambda)$ satisfies conditions $1^\circ, 2^\circ$.
\end{thm}
The smoothness of $q(x)$ is determined by $\omega(\lambda)$ as follows: if $a(x)$ has derivatives of degree $n+4$, then $q(x)$ possesses
$n$ derivatives. Conversely, if $q(x)$ has $n+1$ derivatives, then $a(x)$ has $n+4$ continuous derivatives. Notice also that 
$\int\limits_{-\infty}^0 e^{\sqrt{|\lambda|}x} d\omega(\lambda)$ defines an analytic function on $(0,\infty)$, because its derivatives exists for all $n$:
\[ \int\limits_{-\infty}^0 (\sqrt{|\lambda|})^ne^{\sqrt{|\lambda|}x} d\omega(\lambda) 
\leq \int\limits_{-\infty}^0 e^{\sqrt{|\lambda|}(x+n)} d\omega(\lambda) < \infty,
\]
following from $\sqrt{|\lambda|} \leq e^{\sqrt{|\lambda|}}$.
The existence
of the derivatives of $a(x)$ is directly related to the convergences of the integrals for all $x\in(0,\infty)$ for consecutive values of $n$
\[ \int\limits_1^\infty (\lambda)^n \cos(\sqrt{\lambda}\,x)d\omega(\lambda), \quad
\int\limits_1^\infty \sqrt{\lambda} (\lambda)^n \sin(\sqrt{\lambda}\,x)d\omega(\lambda).
\]

One direction of Theorem \ref{thm:GLMnThm} is classical and very well studied. The second direction of reconstruction of $q(x)$ and $h$ for a given $\rho(\lambda)$
involves the function
\begin{equation} \label{eq:DefF} F(x,y) = \int\limits_{-\infty}^{+\infty} \dfrac{\sin(\sqrt{\lambda}\,x)\sin(\sqrt{\lambda}\,y)}{\lambda} d\omega(\lambda),
\end{equation}
which turns to possess the second mixed derivative
\begin{equation} \label{eq:Deff} f(x,y) = \dfrac{\partial^2 F}{\partial x\partial y}.
\end{equation}
The authors notice that when $\omega(\lambda)$ behaves ``nicely'' at infinity (if for example $var[\omega(\lambda)]<\infty$) $f(x,y)$ can be
represented as follows
\addtocounter{equation}{-1}
\renewcommand{\theequation}{\arabic{equation}'}
\begin{equation} \label{eq:Deff'} f(x,y) = \dfrac{\partial^2 F}{\partial x\partial y} = \int\limits_{-\infty}^{+\infty} \cos(\sqrt{\lambda}\,x)\cos(\sqrt{\lambda}\,y) d\omega(\lambda).
\end{equation}
\renewcommand{\theequation}{\arabic{equation}}
The potential $q(x)$ turns to be
\begin{equation} \label{eq:dKxxGL}
 q(x) = 2 \dfrac{d K(x,x)}{dx},
\end{equation}
where $K(x,y)$, appearing also at \eqref{eq:phicos} satisfies the famous Gelfand-Levitan equation
\begin{equation} \label{eq:GelfandLevitan}
f(x,y) + K(x,y) + \int\limits_{0}^x K(x,t) f(t,y) dt = 0.
\end{equation}

\subsection{Vessels}
In a series of papers of the author \cite{bib:SLVessels, bib:GenVessel, bib:KdVVessels, bib:ENLS} and collaborators
\cite{bib:amv, bib:SchurIEOT} there was developed a similar approach to the inverse scattering 
of many LDEs (like SL, NLS, Canonical systems) using a theory of vessels. In the most general setting in the theory of vessels,
the spectral function $\rho(\lambda)$ is translated into a complex-valued matrix function, belonging
to a special class of functions, parallel to the ``scattering matrix'' $s(\lambda)$ appearing in the work on
Inverse scattering of Faddeyev \cite{bib:FaddeyevII}:
\begin{defn} Let $\sigma_1=\sigma_1^*$ be an invertible $p\times p$ matrix.
Class $\U(\sigma_1)$ consist of $p\times p$ matrix-valued functions $S(\lambda)$ of the complex variable $\lambda$, possessing the following representation:
\begin{equation} \label{eq:S0realized}
S(\lambda) = I - B_0^* \mathbb X_0^{-1} (\lambda I - A)^{-1} B_0 \sigma_1
\end{equation}
where for an auxiliary Hilbert space $\mathcal H$ there are defined operators $B_0:\mathbb C^p\rightarrow\mathcal H$,
$A, \mathbb X_0: \mathcal H\rightarrow\mathcal H$. A general matrix-function $S(\lambda)$,
representable in such a form is called \textbf{realized}. Moreover, the operators are subject to the following assumptions:
\begin{enumerate}
	\item the operator $A$ has a dense domain $D(A)$ with a symmetric spectrum with respect to the imaginary axis.
		$A$ is the generator of a $C_0$ semi-group on $\mathcal H$. Denote the resolvent of A as $R(\lambda)=(\lambda I - A)^{-1}$,
	\item the operator $B_0$ satisfies $R(\lambda) B_0 e \in \mathcal H$ for all $\lambda\not\in\SPEC(A), e\in\mathbb C^p$.
	\item the operator $\mathbb X_0$ is bounded, self adjoint, invertible,
	\item the \textbf{Lyapunov equation} holds for all $\lambda\not\in\SPEC(A)$:
	\begin{equation} \label{eq:LyapunovAt0}
	 \mathbb X_0 R^*(-\bar\lambda) + R(\lambda) \mathbb X_0 - R(\lambda) B_0 \sigma_1 B_0^* R^*(-\bar\lambda) = 0.
	\end{equation}
	In this case $S(\lambda)$ is also called \textbf{symmetric}.
\end{enumerate}
\end{defn}
We call the last equation Lyapunov, because in the case all the operators are bounded, we obtain a usual Lyapunov equation $A \mathbb X + \mathbb X A^* + B_0 \sigma_1 B_0^* = 0$.
It is easy to show using the realization \eqref{eq:S0realized} and the Lyapunov equation only that the symmetry condition
$S^*(-\bar\lambda) \sigma_1 S(\lambda) =  \sigma_1$ holds for all $\lambda\not\in\SPEC(A)$:
\begin{thm}\label{thm:SSym}
Suppose that $S(\lambda)\in\U(\sigma_1$) then $S^*(-\bar\lambda) \sigma_1 S(\lambda) =  \sigma_1$ for all $\lambda\not\in\SPEC(A)$.
\end{thm}
\noindent\textbf{Proof:} Let us plug first the realization formula \eqref{eq:S0realized} for $S(\lambda)$ 
into the symmetry condition:
\begin{equation}\label{eq:SSymProgrs} \begin{array}{lll}
S^*(-\bar\lambda) \sigma_1 S(\lambda) & = (I + \sigma_1B_0^* (\lambda I + A^*)^{-1} \mathbb X_0^{-1} B_0 ) \sigma_1 (I - B_0^* \mathbb X_0^{-1} (\lambda I - A)^{-1} B_0 \sigma_1) \\
& = \sigma_1 + \sigma_1B_0^* (\lambda I + A^*)^{-1} \mathbb X_0^{-1} B_0 \sigma_1 - \sigma_1 B_0^* \mathbb X_0^{-1} (\lambda I - A)^{-1} B_0 \sigma_1 - \\
&\quad \quad -  \sigma_1 B_0^* (\lambda I + A^*)^{-1} \mathbb X_0^{-1} B_0  \sigma_1 B_0^* \mathbb X_0^{-1} (\lambda I - A)^{-1} B_0 \sigma_1
\end{array} \end{equation}
To the last term we apply the Lyapunov equation, rewriting it:
\[ \begin{array}{lll}
\sigma_1 B_0^* (\lambda I + A^*)^{-1} \mathbb X_0^{-1} B_0  \sigma_1 B_0^* \mathbb X_0^{-1} (\lambda I - A)^{-1} B_0 \sigma_1 =  \\
= \sigma_1 B_0^* (\lambda I + A^*)^{-1} \mathbb X_0^{-1} (\lambda I - A)(\lambda I - A)^{-1}B_0  \sigma_1 B_0^* (\lambda I +A^*)^{-1} (\lambda I +A^*)  \mathbb X_0^{-1} (\lambda I - A)^{-1} B_0 \sigma_1 \\
= \sigma_1 B_0^* (\lambda I + A^*)^{-1} \mathbb X_0^{-1} (\lambda I - A)[\mathbb X_0 (\lambda I +A^*)^{-1} 
- (\lambda I - A)^{-1} \mathbb X_0 ] (\lambda I +A^*)  \mathbb X_0^{-1} (\lambda I - A)^{-1} B_0 \sigma_1 \\
= \sigma_1B_0^* (\lambda I + A^*)^{-1} \mathbb X_0^{-1} B_0 \sigma_1 - \sigma_1 B_0^* \mathbb X_0^{-1} (\lambda I - A)^{-1} B_0 \sigma_1
\end{array} \]
If we plug this expression back into \eqref{eq:SSymProgrs} we obtain the symmetry condition. \qed

When $S(\lambda)$ is just analytic at infinity (hence $A$ must be bounded), 
there is a very well known theory of realizations developed in \cite{bib:bgr}. 
For analytic at infinity and symmetric, i.e. satisfying $S^*(-\bar\lambda) \sigma_1 S(\lambda) =  \sigma_1$,
functions there exists a good realization theory using Krein spaces ($\mathcal H$ is a Krein space), 
developed in \cite{bib:KreinReal}\footnote{At the paper \cite{bib:KreinReal} a similar result is proved for functions symmetric with respect to the unit circle,
but it can be translated using Calley transform into $S^*(-\bar\lambda) \sigma_1 S(\lambda) =  \sigma_1$ and was done in \cite{bib:GenVessel, bib:SchurIEOT}}.
Such a realization is then translated into a function in $\U(\sigma_1)$.
In order to acquit the name of ``scattering data'' for $S(\lambda)$ we show how to construct the ``potential''
(analogously to SL case), which we call $\gamma_*(x)$. For this we will need to fix two additional $p\times p$ matrices
$\sigma_2=\sigma_2^*$ and $\gamma=-\gamma^*$:
\begin{enumerate}
	\item Suppose that there exists a solution $B(x)$ satisfying $B(0)=B_0$ and the equation
	\begin{equation} \label{eq:DB}
		0  =  \frac{d}{dx} (R(\lambda) B(x)\sigma_1 e ) + 
		A R(\lambda) B(x)\sigma_2 e  + R(\lambda) B(x) \gamma e, \quad 
		\forall \lambda\not\in\SPEC(A), e\in\mathbb C^p.
	\end{equation}
	This equation is solvable because the coefficients $\sigma_1, \sigma_2, \gamma$ are constant and $A$ is
	a generator of a $C_0$ semigroup. We will solve it in an exact manner for the SL case. For this equation to 
	hold we also demand that the following regularity assumptions hold:
	\begin{eqnarray}
	\label{eq:ResBsigma2} \forall\lambda\not\in\SPEC(A): R(\lambda) B(x)\sigma_2\mathbb C^p\subseteq D(A), \\ 
	\label{eq:ResBgamma} \forall\lambda\not\in\SPEC(A): R(\lambda) B(x)\gamma\mathbb C^p\subseteq\mathcal H.
	\end{eqnarray}
\item Solve for $\mathbb X(x)$, satisfying $\mathbb X(0)=\mathbb X_0$ and:
\begin{equation}
\label{eq:DX} \frac{d}{dx} \mathbb X(x)  =  B(x) \sigma_2 B^*(x).
\end{equation}
Notice that $B(x)\sigma_2\mathbb C^p\subseteq\mathcal H$ follows from \eqref{eq:ResBsigma2}.
\item Define $\gamma_*(x)$ on $\mathrm I$ by
\begin{equation}
\label{eq:Linkage}
\gamma_*(x)  =  \gamma + \sigma_2 B^*(x) \mathbb X^{-1}(x) B(x) \sigma_1 
 - \sigma_1 B^*(x) \mathbb X^{-1}(x) B(x) \sigma_2.
\end{equation}
\end{enumerate}
The main reason, why we call $S(\lambda)$ as the ``scattering data'' is the fact that $\gamma_*(x)$ (generalized
potential) is uniquely determined from $S(\lambda)$ by this construction. 
With some additional requirements, there is also a unique $S(\lambda)$ for each
$\gamma_*(x)$. Still, the question of reconstructing $S(\lambda)$ from the given $\gamma_*(x)$ is the fundamental 
problem of the theory of vessels, solution of which we are going to present here for the Sturm Liouville case under a certain regularity
assumption (Theorem \ref{thm:VesselFromq}).  Moreover, the function 
\begin{equation} \label{eq:SlxRealized}
 S(\lambda,x) = I - B^*(x)\mathbb X^{-1}(x) (\lambda I - A)^{-1} B(x) \sigma_1
\end{equation} 
acts as a B\" acklund transformation, because (see Theorem \ref{thm:Backlund})
multiplication by it maps solutions $u(\lambda,x)$ of the \textbf{input} 	LDE with the spectral parameter $\lambda$:
\begin{equation} \label{eq:InCC}
		-\sigma_1\dfrac{\partial}{\partial x}u(\lambda,x) + (\sigma_2 \lambda + \gamma)u(\lambda,x) = 0
\end{equation}
to solutions $y(\lambda,x)=S(\lambda,x)u(\lambda,x)$ of the \textbf{output}	LDE with the same spectral parameter:
\begin{equation} \label{eq:OutCC}
		-\sigma_1\dfrac{\partial}{\partial x} y(\lambda,x) + (\sigma_2 \lambda + \gamma_*(x))y(\lambda,x) = 0.
\end{equation}
Notice also that we relate in this manner solutions of an LDE with constant coefficients to solutions of an LDE with
variable coefficients analogously to the scattering theory.

Such functions are called \textit{transfer functions} of vessels, defined in the next
\begin{defn} \label{def:Vessel}The collection of operators and spaces 
\begin{equation} \label{eq:DefV}
\mathfrak{V} = (A, B(x), \mathbb X(x); \sigma_1, \sigma_2, \gamma, \gamma_*(x);
\mathcal{H},\mathbb C^p; \mathrm I),
\end{equation}
is called a \textbf{vessel}, if $\mathbb X(x), A:\mathcal H\rightarrow\mathcal H$, 
$B(x):\mathbb C^p\rightarrow\mathcal H$ are differentiable linear operators, subject to the following conditions:
\begin{enumerate}
	\item the operator $A$ has a dense domain $D(A)$ with a symmetric spectrum with respect to the imaginary axis.
		$A$ is the generator of a $C_0$ semi-group on $\mathcal H$,
	\item $B(x)$ satisfies regularity assumptions \eqref{eq:ResBsigma2}, \eqref{eq:ResBgamma}, and the equation \eqref{eq:DB}.
	\item $\mathbb X(x)$ is bounded, self adjoint, invertible on $\mathrm I$ and satisfies \eqref{eq:DX},
	\item the \textbf{Lyapunov equation} holds for all $x\in\mathrm I, \lambda\not\in\SPEC(A)$:
	\begin{equation} \label{eq:Lyapunov}
		 \mathbb X(x) R^*(-\bar\lambda) + R(\lambda) \mathbb X(x) - R(\lambda) B(x) \sigma_1 B(x)^*R^*(-\bar\lambda) = 0.
	\end{equation}
	\item $\gamma_*(x)$ satisfies \eqref{eq:Linkage}.
\end{enumerate}
\end{defn}
The class of transfer functions is defined as follows
\begin{defn} Class $\I(\sigma_1)$ consist of $p\times p$ matrix-valued (transfer) functions $S(\lambda,x)$ of the complex variable $\lambda$ and $x\in\mathrm I$ for an interval $\mathrm I=[a,b]$ possessing the following representation:
\begin{equation} \label{eq:Srealized}
 S(\lambda,x) = I - B^*(x) \mathbb X^{-1}(x) (\lambda I - A)^{-1} B(x) \sigma_1 
\end{equation}
where for an auxiliary Hilbert space $\mathcal H$, the operators $B(x):\mathbb C^p\rightarrow\mathcal H$, $\mathbb X(x), A: \mathcal H\rightarrow\mathcal H$ constitute a vessel
for some $\sigma_2,\gamma_*(x)$ on $\mathrm I$.
\end{defn}
we want to prove the B\" acklund transformation theorem, presented earlier. Before we do that few techincal results are needed:
\begin{lemma}[\textbf{permanence of the Lyapunov equation}] \label{lemma:Redund}
Suppose that $B(x)$ satisfies \eqref{eq:DB} and $\mathbb X(x)$ satisfies \eqref{eq:DX}, then if the Lyapunov equation \eqref{eq:Lyapunov}
holds for a fixed $x_0$, then it holds for all $x$. If $\mathbb X(x_0) = \mathbb X^*(x_0)$ then
$\mathbb X(x)$ is self-adjoint for all $x$.
\end{lemma}
\noindent\textbf{Proof:} Let us differentiate Lyapunov equation \eqref{eq:Lyapunov}:
\[ \begin{array}{lll}
\dfrac{d}{dx} [-\mathbb X(x)(\lambda I + A^*)^{-1} + (\lambda I - A)^{-1}\mathbb X(x) + (\lambda I - A)^{-1} B(x) \sigma_1 B^*(x) (\lambda I + A^*)^{-1}] = \\
= - B(x)\sigma_2B^*(x)(\lambda I + A^*)^{-1} + (\lambda I - A)^{-1}B(x)\sigma_2B^*(x)  + \\
\quad \quad + \dfrac{d}{dx} [(\lambda I - A)^{-1} B(x) \sigma_1 B^*(x) (\lambda I + A^*)^{-1}] = \\
= - B(x)\sigma_2B^*(x)(\lambda I + A^*)^{-1} + (\lambda I - A)^{-1}B(x)\sigma_2B^*(x)  + \\
\quad \quad - A (\lambda I - A)^{-1} B(x) \sigma_2 B^*(x) (\lambda I + A^*)^{-1} - (\lambda I - A)^{-1} B(x)\cancel{\gamma} B^*(x) (\lambda I + A^*)^{-1} - \\
\quad \quad \quad \quad  - (\lambda I - A)^{-1} B(x) \sigma_2 B^*(x) (\lambda I + A^*)^{-1} A^* -
 (\lambda I - A)^{-1} B(x) \cancel{\gamma^*} B^*(x) (\lambda I + A^*)^{-1}
\end{array} \] 
The terms involving $\gamma$ are canceled, because $\gamma+\gamma^*=0$, by the assumption on it. So this becomes, after adding and subtracting $\lambda I$ next to $A$ and to $A^*$:
\[ \begin{array}{lll}
= - B(x)\sigma_2B^*(x)(\lambda I + A^*)^{-1} + (\lambda I - A)^{-1}B(x)\sigma_2B^*(x) + \\
\quad \quad - (A-\lambda I + \lambda I) (\lambda I - A)^{-1} B(x) \sigma_2 B^*(x) (\lambda I + A^*)^{-1} - \\
\quad \quad \quad \quad  - (\lambda I - A)^{-1} B(x) \sigma_2 B^*(x) (\lambda I + A^*)^{-1} (A^* +\lambda I  - \lambda I) = \\
= - B(x)\sigma_2B^*(x)(\lambda I + A^*)^{-1} + (\lambda I - A)^{-1}B(x)\sigma_2B^*(x) + \\
\quad \quad + B(x) \sigma_2 B^*(x) (\lambda I + A^*)^{-1} - \cancel{\lambda} (\lambda I - A)^{-1} B(x) \sigma_2 B^*(x) (\lambda I + A^*)^{-1} + \\
\quad \quad \quad \quad  - (\lambda I - A)^{-1} B(x) \sigma_2 B^*(x)  + (\lambda I - A)^{-1} B(x) \sigma_2 B^*(x) (\lambda I + A^*)^{-1} \cancel{\lambda} = \\
= 0.
\end{array} \] 
So, it is a constant and the result follows. 

Since the derivative \eqref{eq:DX} $\dfrac{d}{dx}\mathbb X(x) = B(x) \sigma_2(x) B^*(x)$
is self-adjoint, $\mathbb X(x)$ will be self-adjoint, once $\mathbb X(x_0)$ is. \qed
\begin{lemma} \label{lemma:X-1B}
Suppose that $B(x)$ satisfies \eqref{eq:DB} and $\mathbb X(x)$ satisfies \eqref{eq:DX}, then
\begin{multline} \label{eq:X-1B}
\dfrac{d}{dx}[(\lambda I + A^*)^{-1}\mathbb X^{-1}(x) B(x) \sigma_1] =  \\
-\lambda (\lambda I + A^*)^{-1}\mathbb X^{-1}(x) B(x) \sigma_2 + \mathbb X^{-1}(x) B(x) \sigma_2 - (\lambda I + A^*)^{-1}\mathbb X^{-1}(x) B(x) \gamma_*(x).
\end{multline}
\end{lemma}
\noindent\textbf{Remark: 1. } Notice that if $A$ is bounded it is equivalent to
\begin{equation} \label{eq:DX-1BBounded} \dfrac{d}{dx}[\mathbb X^{-1} B(x) \sigma_1] = 
A^* \mathbb X^{-1} B(x) \sigma_2 - \mathbb X^{-1} B(x) \gamma_*(x)
\end{equation}
\textbf{2.} Substituting $\lambda$ with $-\bar\lambda$ in \eqref{eq:X-1B} and taking adjoint, we will obtain that the following
equality holds
\begin{multline} \label{eq:DB*X-1}
\dfrac{d}{dx}[\sigma_1B^*(x)\mathbb X^{-1}(x) R(\lambda)] =  \\
 \lambda \sigma_2 B^*(x)\mathbb X^{-1}(x) R(\lambda) - \sigma_2 B^*(x)\mathbb X^{-1}(x)-
\gamma_*(x) B^*(x)\mathbb X^{-1}(x)  R(\lambda).
\end{multline}
\noindent\textbf{Proof of Lemma \ref{lemma:X-1B}:} Let us differentiate $(\lambda I + A^*)^{-1}\mathbb X^{-1} B(x) \sigma_1$ using \eqref{eq:DB}, \eqref{eq:DX}:
\[ \begin{array}{lll}
\dfrac{d}{dx}[(\lambda I + A^*)^{-1}\mathbb X^{-1} B(x) \sigma_1] = \\
= (\lambda I + A^*)^{-1}\dfrac{d}{dx}[\mathbb X^{-1}] B(x) \sigma_1 + 
(\lambda I + A^*)^{-1}\mathbb X^{-1} (\lambda I -A) \dfrac{d}{dx}[(\lambda I -A)^{-1} B(x) \sigma_1] = \\
= -(\lambda I + A^*)^{-1}\mathbb X^{-1} B(x)\sigma_2B^*(x) \mathbb X^{-1} B(x) \sigma_1 +\\
\quad \quad + (\lambda I + A^*)^{-1}\mathbb X^{-1}(\lambda I -A) [- A (\lambda I -A)^{-1} B(x) \sigma_2 - (\lambda I -A)^{-1} B(x) \gamma] .
\end{array} \]
Notice that $\lambda I - A$ is a bounded operator for $\lambda\not\in\SPEC(A)$. Then rearranging
\[ \begin{array}{lll}
= -(\lambda I + A^*)^{-1}\mathbb X^{-1} B(x)\sigma_2B^*(x) \mathbb X^{-1} B(x) \sigma_1 - (\lambda I + A^*)^{-1}\mathbb X^{-1}B(x) \gamma -\\
\quad \quad - (\lambda I + A^*)^{-1}\mathbb X^{-1}(\lambda I -A) A (\lambda I -A)^{-1} B(x) \sigma_2.
\end{array} \]
If we add and subtract $\lambda I$ next to $A$ at the last term, we obtain
\begin{equation} \label{eq:DX-1B*mid} \begin{array}{lll}
= -(\lambda I + A^*)^{-1}\mathbb X^{-1} B(x)\sigma_2B^*(x) \mathbb X^{-1} B(x) \sigma_1 - (\lambda I + A^*)^{-1}\mathbb X^{-1}B(x) \gamma -\\
\quad \quad - (\lambda I + A^*)^{-1}\mathbb X^{-1}(\lambda I -A) (A -\lambda I + \lambda I )(\lambda I -A)^{-1} B(x) \sigma_2 \\
= -(\lambda I + A^*)^{-1}\mathbb X^{-1} B(x) [\gamma+\sigma_2B^*(x) \mathbb X^{-1} B(x) \sigma_1 ] - \\
\quad \quad - \lambda (\lambda I + A^*)^{-1}\mathbb X^{-1} B(x) \sigma_2 + (\lambda I + A^*)^{-1}\mathbb X^{-1}(\lambda I -A) B(x) \sigma_2 \\
\end{array} \end{equation}
A simple consequence of the Lyapunov equation \eqref{eq:Lyapunov} is as follows
\[ (\lambda I + A^*)^{-1}\mathbb X^{-1}(\lambda I -A) = \mathbb X^{-1}(x) + (\lambda I + A^*)^{-1} \mathbb X^{-1}(x) B(x)\sigma_1B(x) \mathbb X^{-1}(x)
\]
So if we plug this expression into the last term of \eqref{eq:DX-1B*mid}, we obtain rearranging
\[ \begin{array}{lll}
\dfrac{d}{dx}[(\lambda I + A^*)^{-1}\mathbb X^{-1} B(x) \sigma_1] = \\
=  -(\lambda I + A^*)^{-1}\mathbb X^{-1} B(x) [\gamma+\sigma_2B^*(x) \mathbb X^{-1} B(x) \sigma_1 ] - \\
\quad \quad - \lambda (\lambda I + A^*)^{-1}\mathbb X^{-1} B(x) \sigma_2 + [\mathbb X^{-1}(x) + (\lambda I + A^*)^{-1} \mathbb X^{-1}(x) B(x)\sigma_1B(x) \mathbb X^{-1}(x)] B(x) \sigma_2 =\\
=  -(\lambda I + A^*)^{-1}\mathbb X^{-1} B(x) [\gamma+\sigma_2B^*(x) \mathbb X^{-1} B(x) \sigma_1 - \sigma_1B(x) \mathbb X^{-1}(x) B(x) \sigma_2] - \\
\quad \quad - \lambda (\lambda I + A^*)^{-1}\mathbb X^{-1} B(x) \sigma_2 + \mathbb X^{-1}(x)B(x) \sigma_2 \\
= -(\lambda I + A^*)^{-1}\mathbb X^{-1} B(x) \gamma_*(x) - \\
\quad \quad - \lambda (\lambda I + A^*)^{-1}\mathbb X^{-1} B(x) \sigma_2 + \mathbb X^{-1}(x)B(x) \sigma_2,
\end{array} \]
using \eqref{eq:Linkage} at the final step.
\qed

Now we have all the ingredients of the following Theorem. This theorem has its origins at the work of M. Liv\c sic \cite{bib:Vortices} and was
proved for bounded operators in \cite{bib:GenVessel, bib:ENLS, bib:SchurIEOT}. Now we present a generalization of these results for the case of unbounded
operator $A$.
\begin{thm}[Vessel=B\" acklund transformation] \label{thm:Backlund}
Let $\mathfrak{V}$ be a vessel defined in \eqref{eq:DefV} and satisfying all the conditions of Definition \ref{def:Vessel}. Fix
$\lambda\not\in\SPEC(A)$ and let $u(\lambda,x)$ be a solution of the input LDE \eqref{eq:InCC}. Then the function $y(\lambda,x)=S(\lambda,x)u(\lambda,x)$
is differentiable and satisfies the output LDE \eqref{eq:OutCC}.
\end{thm}
\noindent\textbf{Proof:} Let us fix $\lambda\not\in\SPEC(A)$ and a solution $u(\lambda,x)$ of \eqref{eq:InCC}.
Then for $y(\lambda,x) = S(\lambda,x) u(\lambda,x)$ we calculate:
\[ \begin{array}{lll}
\sigma_1 \dfrac{d}{dx} y(\lambda,x)  & = \sigma_1 \dfrac{d}{dx} [(I - B^*(x)\mathbb X^{-1}(x) (\lambda I - A)^{-1} B(x)\sigma_1) u(\lambda,x)] =\\
& = \sigma_1 \dfrac{d}{dx} u(\lambda,x) - \sigma_1 \dfrac{d}{dx} [ B^*(x)\mathbb X^{-1}(x) (\lambda I - A)^{-1} B(x)\sigma_1 u(\lambda,x)] \\
& = (\sigma_2\lambda+\gamma)u(\lambda,x) - \sigma_1 \dfrac{d}{dx} [ B^*(x)\mathbb X^{-1}(x)] \, (\lambda I - A)^{-1} B(x)\sigma_1 u(\lambda,x) \\
& \quad \quad - \sigma_1B^*(x)\mathbb X^{-1}(x) \dfrac{d}{dx} [(\lambda I - A)^{-1} B(x)] u(\lambda,x)  \\
&  \quad \quad \quad \quad - \sigma_1B^*(x)\mathbb X^{-1}(x) (\lambda I - A)^{-1} B(x) \sigma_1 \dfrac{d}{dx} u(\lambda,x) .
\end{array} \]
Using \eqref{eq:DB*X-1} and \eqref{eq:DB} it becomes
\[ \begin{array}{llll}
\sigma_1 \dfrac{d}{dx} y(\lambda,x)   =
(\sigma_2\lambda+\gamma)u(\lambda,x) - \\
-  [\lambda \sigma_2 B^*(x)\mathbb X^{-1}(x)(\lambda I - A)^{-1} - \sigma_2 B^*(x)\mathbb X^{-1}(x)-
\gamma_*(x) B^*(x)\mathbb X^{-1}(x) (\lambda I - A)^{-1}] B(x)\sigma_1 u(\lambda,x) +\\
\quad \quad + \sigma_1B^*(x)\mathbb X^{-1}(x) [A (\lambda I - A)^{-1} B(x)\sigma_2 + (\lambda I - A)^{-1} B(x)\cancel{\gamma}] u(\lambda,x) \\
\quad \quad \quad \quad -\sigma_1B^*(x)\mathbb X^{-1}(x) (\lambda I - A)^{-1} B(x) (\sigma_2\lambda+\cancel\gamma) u(\lambda,x) =
\end{array} \]
Let us combine the last two terms and rearrange the rest:
\[ \begin{array}{llll}
= (\sigma_2\lambda + \gamma)u(\lambda,x) - \\
-  \lambda \sigma_2 B^*(x)\mathbb X^{-1}(x)(\lambda I - A)^{-1}  B(x)\sigma_1 u(\lambda,x) 
- \sigma_2 B^*(x)\mathbb X^{-1}(x)B(x)\sigma_1 u(\lambda,x) - \\
\gamma_*(x) B^*(x)\mathbb X^{-1}(x) (\lambda I - A)^{-1}] B(x) \sigma_1 u(\lambda,x) + \\
\quad \quad + \sigma_1B^*(x)\mathbb X^{-1}(x) (A -\lambda I)(\lambda I - A)^{-1} B(x)\sigma_2 u(\lambda,x) \\
= \sigma_2\lambda (I - B^*(x)\mathbb X^{-1}(x)(\lambda I - A)^{-1}  B(x)\sigma_1) u(\lambda,x) + \\
+ [\gamma + \sigma_2 B^*(x)\mathbb X^{-1}(x)B(x)\sigma_1- \sigma_1B^*(x)\mathbb X^{-1}(x) B(x)\sigma_2]  u(\lambda,x)- \\
- \gamma_*(x) B^*(x)\mathbb X^{-1}(x) (\lambda I - A)^{-1} B(x) \sigma_1 u(\lambda,x)
\end{array} \]
Using \eqref{eq:Linkage} and the definition of $S(\lambda,x)$ we obtain that it is
\[ \begin{array}{llll}
= \sigma_2\lambda S(\lambda,x)u(\lambda,x) + \gamma_*(x) u(\lambda,x) - \\
- \gamma_*(x) B^*(x)\mathbb X^{-1}(x) (\lambda I - A)^{-1} B(x) \sigma_1 u(\lambda,x) = \\
= \sigma_2\lambda S(\lambda,x)u(\lambda,x) + \gamma_*(x) S(\lambda,x) u(\lambda,x) =\\
= (\sigma_2\lambda + \gamma_*(x)) S(\lambda,x) u(\lambda,x) = \\
= (\sigma_2\lambda + \gamma_*(x)) y(\lambda,x).
\end{array} \]
\qed

We have seen that there is a ``simple'' construction, called \textit{standard}, of functions in $\I(\sigma_1)$,
starting from a function in $\U(\sigma_1)$ by solving differential equations. Notice also that in the case $\mathbb X_0=I$ (the identity operator)
we have a very well developed theory \cite{bib:BL} of so called ``characteristic functions'' of the operator $A$, which has a finite dimensional real part:
$A + A^* = -B\sigma_1B^*$. But in general, when $\mathbb X_0$ is not identity, it is only known that $A$ has a spectrum, symmetric with respect to the imaginary
axis. Existence of the vessel and its transfer function relies on the invertability of the operator $\mathbb X(x)$. In order to investigate the existence of the inverse
for $\mathbb X(x)$ notice that from \eqref{eq:DX}
\[ \mathbb X(x) = \mathbb X_0 + \int\limits_0^x B(y) \sigma_2 B^*(y) dy
\]
it follows that
\[ \mathbb X_0^{-1} \mathbb X(x) = I + \mathbb X^{-1}_0 \int\limits_{x_0}^x B(y) \sigma_2 B^*(y) dy.
\]
Since $\sigma_2$ has finite rank for $\dim \mathcal E<\infty$, 
this expression is of the form $I + T$, for a trace class operator $T$ and since 
$\mathbb X_0$ is an invertible operator, there exists a non trivial interval (of length at least $\dfrac{1}{\|\mathbb X_0^{-1}\|}$) on which $\mathbb X(x)$ and $\tau(x)$ are defined. Recall \cite{bib:GKintro} that a function $F(x)$ from (a, b) into the group G (the set of bounded invertible operators on H of the form I + T, for
a trace-class operator $T$) is said to be differentiable if $F(x) -I$ is \textit{differentiable} as a map into the trace-class operators. In our case,
\[ \dfrac{d}{dx} (\mathbb X_0^{-1}\mathbb X(x)) = 
\mathbb X_0^{-1} \dfrac{d}{dx} \mathbb X(x) =
\mathbb X_0^{-1} B(x)\sigma_2B^*(x)
\]
exists in trace-class norm. This leads us to the following
\begin{defn} \label{def:Tau} For a given vessel $\mathfrak{V}$ \eqref{eq:DefV} the tau function $\tau(x)$ is defined as
\begin{equation} \label{eq:Tau} \tau = \det (\mathbb X_0^{-1} \mathbb X(x)).
\end{equation}
The vessel exists at a neighborhood of every $x$ for which $\tau(x)\neq 0$.
\end{defn}

Israel Gohberg and Mark Krein \cite[formula 1.14 on p. 163]{bib:GKintro}
proved that if $\mathbb X_0^{-1}\mathbb X(x)$ is a differentiable function
into G, then $\tau(x) = \SP (\mathbb X_0^{-1}\mathbb X(x))$
\footnote{$\SP$ - stands for the trace in the infinite dimensional space.} is a differentiable map into $\mathbb C^*$ with
\begin{multline} \label{eq:GKform}
\dfrac{\tau'}{\tau}  = \SP (\big(\mathbb X_0^{-1} \mathbb X(x)\big)^{-1} 
\dfrac{d}{dx} \big(\mathbb X_0^{-1} \mathbb X(x)\big)) = \SP (\mathbb X(x)' \mathbb X^{-1}(x)) = \\
= \SP (B(x)\sigma_2 B^*(x) \mathbb X^{-1}(x)) =
\TR (\sigma_2 B^*(x) \mathbb X^{-1}(x)B(x)).
\end{multline}
Differentiating this expression, we obtain that
\[ (\dfrac{\tau'}{\tau})' = \dfrac{\tau''}{\tau} - (\dfrac{\tau'}{\tau})^2 =
\dfrac{d}{dx} \TR (\sigma_2 B^*(x) \mathbb X^{-1}(x)B(x)).
\]
Using vessel conditions, since $B(x)$, $\mathbb X^{-1}(x)$ are differentiable bounded operators in the case
$A\mathbb X^{-1}(x)B(x)$ exists, or in the case it is canceled (SL case) we obtain that
\begin{equation} \label{eq:tau'tau'}
 \begin{array}{llll}
\dfrac{d}{dx} \TR (\sigma_2 B^*(x) \mathbb X^{-1}(x)B(x)) = \\
=  \TR (\sigma_2 \sigma_1^{-1}(-\sigma_2 B^*(x)A^* - \gamma^* B^*(x))\mathbb X^{-1}(x)B(x)) - \\
\quad - \TR(\sigma_2 B^*(x) \mathbb X^{-1}(x)B(x) \sigma_2 B^*(x) \mathbb X^{-1}(x) B(x))
+ \TR(\sigma_2 B^*(x) \mathbb X^{-1}(x) (-A B(x)\sigma_2 - B(x)\gamma) \sigma_1^{-1}) = \\
= \TR (\sigma_2 \sigma_1^{-1} \sigma_2 B^*(x)[ -A^* \mathbb X^{-1} - \mathbb X^{-1} A]B(x)) - 
\TR( [\sigma_2 \sigma_1^{-1}\gamma^* + \gamma\sigma_1^{-1}\sigma_2] B^*(x)\mathbb X^{-1}(x)B(x)) - \\
\quad - \TR(\sigma_2 B^*(x) \mathbb X^{-1}(x)B(x) \sigma_2 B^*(x) \mathbb X^{-1}(x) B(x))= \\
= \TR (\sigma_2 \sigma_1^{-1} \sigma_2 B^*(x) \mathbb X^{-1}B(x)\sigma_1 B^*(x) \mathbb X^{-1} B(x)) - 
\TR( [\sigma_2 \sigma_1^{-1}\gamma^* + \gamma\sigma_1^{-1}\sigma_2] B^*(x)\mathbb X^{-1}(x)B(x)) - \\
\quad - \TR(\sigma_2 B^*(x) \mathbb X^{-1}(x)B(x) \sigma_2 B^*(x) \mathbb X^{-1}(x) B(x)).
\end{array} \end{equation}
In the special case of SL vessel parameters, we obtain that equations \eqref{eq:InCC}, \eqref{eq:OutCC} are equivalent to \eqref{eq:SL}. Let us explain it
in more details.
\begin{defn} \label{def:SLparam}
The Sturm Liouville (SL) vessel parameters are defined as follows
\[ \sigma_1 = \bbmatrix{0 & 1 \\ 1 & 0}, \quad
\sigma_2 = \bbmatrix{1 & 0 \\ 0 & 0}, \quad
\gamma = \bbmatrix{0 & 0 \\ 0 & i}.
\]
\end{defn}
Suppose that we are given a SL vessel $\mathfrak{V}$, in other words, $\mathfrak{V}$ \eqref{eq:DefV} is defined
for the SL vessel parameters. Denote a differentiable $2\times 2$ matrix function 
$H_0 = B^*(x) \mathbb X^{-1}(x) B(x) = \bbmatrix{a&b\\b^*&d}$. 
Then from \eqref{eq:GKform} it follows that
$\dfrac{\tau'}{\tau}  = \TR (\sigma_2 B^*(x) \mathbb X^{-1}(x)B(x)) = a$ is the $1,1$ entry of $H_0$.
Using simple calculations it follows that
\[\sigma_2\sigma_1^{-1}\sigma_2=0, \quad \TR (\sigma_2 H_0\sigma_2H_0) = (\TR(\sigma_2H_0))^2, \quad
\sigma_2 \sigma_1^{-1}\gamma^* + \gamma\sigma_1^{-1}\sigma_2 = \bbmatrix{0&-i\\i&0} , \] 
and we obtain from \eqref{eq:tau'tau'} that
\[ \dfrac{\tau''}{\tau} = - \TR(\bbmatrix{0&-i\\i&0} H_0) = i (b^*-b).
\]
Notice that the terms involving operator $A$ are canceled. Moreover, we obtain that 
\[ \gamma_*(x) = \gamma + \sigma_2 H_0 \sigma_1 - \sigma_1 H_0 \sigma_2 =
\bbmatrix{b-b^* & a \\ -a & i} = \bbmatrix{i \dfrac{\tau''}{\tau} & \dfrac{\tau'}{\tau} \\ -\dfrac{\tau'}{\tau} & i}.
\]
Thus we obtain the following lemma (appearing already in \cite[Proposition 3.2 ]{bib:SLVessels})
\begin{lemma} For SL vessel parameters, the following formula for $\gamma_*(x)$ holds
\[ \gamma_*(x) = \bbmatrix{i \dfrac{\tau''}{\tau} & \dfrac{\tau'}{\tau} \\ -\dfrac{\tau'}{\tau} & i}.
\]
\end{lemma}

Analogously to \cite[Section 3.1.1]{bib:SLVessels}, simple calculations show that denoting $u(\lambda,x) = \bbmatrix{u_1(\lambda,x)\\ u_2(\lambda,x)}$ we shall obtain that the input compatibility
condition \eqref{eq:InCC} is equivalent to
\[ \left\{ \begin{array}{lll}
-\frac{\partial^2}{\partial x^2} u_1(\lambda,x) = -i\lambda u_1(\lambda,x), \\
u_2(\lambda,x) = - i \frac{\partial}{\partial x} u_1(\lambda,x).
\end{array}\right.
\]
The output $y(\lambda,x) = \bbmatrix{y_1(\lambda,x)\\y_2(\lambda,x)} = S(\lambda,x) u(\lambda,x)$ satisfies the output
equation \eqref{eq:OutCC}, which is equivalent to
\[ \left\{ \begin{array}{lll}
-\frac{\partial^2}{\partial x^2} y_1(\lambda,x) - 2 \dfrac{d^2}{dx^2} [\ln(\tau(x))]  y_1(\lambda,x) = -i\lambda y_1(\lambda,x), \\
y_2(\lambda,x) = - i [ \frac{\partial}{\partial x} y_1(\lambda,x) + \dfrac{\tau'}{\tau}y_1(\lambda,x)].
\end{array}\right.
\]
Observing the first coordinates $u_1(\lambda,x), y_1(\lambda,x)$ of the vector-functions $u(\lambda,x),y(\lambda,x)$ we can see that multiplication by
$S(\lambda,x)$ maps solution of the trivial SL equation (i.e. $q(x)=0$) to solutions of the more complicated one,
defined by the potential 
\begin{equation} \label{eq:qbeta}
q(x) = - 2  \dfrac{d^2}{dx^2} [\ln(\tau(x))],
\end{equation}
which can be considered as an analogue of the scattering theory.

\section{Construction of a vessel using Gelfand-Levitan theory}
Throughout of this section we suppose that $h=0$ at the initial value problem \eqref{eq:SLinit}. This will enable us to use the
formulas \eqref{eq:phicos} and \eqref{eq:cosphi}.
Suppose that we are given a potential $q(x)$ possessing a spectral function $\rho(\mu)$, for which
$f(x,y)$ satisfies the condition \eqref{eq:Deff'}:
\[ f(x,y) = \dfrac{\partial^2 F}{\partial x\partial y} = \int\limits_{-\infty}^{+\infty} \cos(\sqrt{\mu}\,x)\cos(\sqrt{\mu}\,y) d\omega(\mu).
\]
From here it follows (for $x=y$) that $\cos(\sqrt{\mu}\,x)$ is in the space
\[
 \mathcal H = L^2_\omega = \{ g(\mu) \mid \int\limits_{\mathbb R} |g(\mu)|^2 d\omega(\mu) < \infty\}, \\
\]
and the functional $[\cos(\sqrt{\lambda}\,x)]^*$, defined by this function:
\[ [\cos(\sqrt{\mu}\,x)]^* g(\mu) = \int\limits_{-\infty}^{+\infty} \cos(\sqrt{\mu}\,x) g(\mu) d\omega(\mu)
\]
is well defined on $\mathcal H$. Moreover, by formula \eqref{eq:phicos} we also obtain that since the kernel $K(x,t)$ is continuous,
the functional $[\phi(x,\lambda)]^*$ is well defined.

We are going to present a construction of a vessel $\mathfrak V$ realizing the potential $q(x)$.
Let us define the operator $A = i \lambda$ as the multiplication operator on $i\lambda$. It is usually an unbounded
operator with the obvious domain. Define next\footnote{for negative $\lambda$ we understand $\cos(\sqrt{\lambda})$ as
$\cosh(\sqrt{|\lambda|})$ and similarly for the sinus}
\begin{equation} \label{eq:DefBdomega}
 B(x) = \bbmatrix{\cos(\sqrt{\mu} x) & i \sqrt{\mu}\sin(\sqrt{\mu}x)}
\end{equation}
Then for $\lambda\not\in\SPEC(A)$, $(\lambda I - A)^{-1} B(x) = \dfrac{1}{\lambda-i\mu}\bbmatrix{\cos(\sqrt{\mu} x) & i \sqrt{\mu}\sin(\sqrt{\mu}x)}$. 
It is almost readable that the conditions 
\eqref{eq:ResBsigma2}, \eqref{eq:ResBgamma} are fulfilled and the equation \eqref{eq:DB} hold:
\begin{multline*}
\dfrac{d}{dx} (\lambda I - A)^{-1} B(x)\sigma_1 =  (\lambda I - A)^{-1} \bbmatrix{ i \lambda \cos(\sqrt{\lambda}x) & -\sqrt{\lambda}\sin(\sqrt{\lambda} x)}
 = \\
 = - A(\lambda I - A)^{-1} B(x)\sigma_2 - (\lambda I - A)^{-1}B(x)\gamma.
\end{multline*}
As a result, we obtain that
\begin{equation} \label{eq:OmeagaDef}
f(x,y) = \bbmatrix{1&0}B^*(x) B(y) \bbmatrix{1\\0} = \bbmatrix{1&0}B^*(x) \mathbb X_0^{-1} B(y) \bbmatrix{1\\0}.
\end{equation}
Define next operator $\mathbb X(x):\mathcal H\rightarrow\mathcal H$ using an arbitrary $g\in\mathcal H$. Here we present how to calculate 
the resulting function $\mathbb X(x) g$ at point $\lambda$:
\begin{equation} \label{eq:XDefInt}
 (\mathbb X(x) g)(\lambda) = [I + \int\limits_0^x B(y)\sigma_2B^*(y) dy ]g =
g(\lambda) +  \int\limits_0^x \cos(\sqrt{\lambda}y) \int\limits_{\mathbb R} \cos(\sqrt{\mu}y) g(\mu) d\omega(\mu)
\end{equation}
notice that this operator is of the form $I + P$ for a positive, trace-class operator $P$,
so it is invertible on $[0,\infty)$.
From these formula we obtain an explicit representation for $K(x,y)$:
\begin{lemma} The following equality holds:
\begin{equation} \label{eq:KxtBasic}
K(x,y) = - \int\limits_{-\infty}^{+\infty} \phi(x,\lambda) \cos(\sqrt{\lambda}\,y) d\omega(\lambda).
\end{equation}
where $\phi(x,\lambda)$ is the solution of \eqref{eq:SL} satisfying $\phi(0,\lambda)=1$, $\phi'(0,\lambda)=0$ \eqref{eq:SLinit}.
\end{lemma}
\noindent\textbf{Proof:} Plugging the expression \eqref{eq:phicos} for $\phi(x,\lambda)$ we obtain
\[ \begin{array}{llll}
\int\limits_{-\infty}^{+\infty} \phi(x,\lambda) \cos(\sqrt{\lambda}\,y) d\omega(\lambda) =
\int\limits_{-\infty}^{+\infty} [\cos(\sqrt{\lambda}\,x) + \int\limits_0^x K(x,t) \cos(\sqrt{\lambda}\,t) dt] \cos(\sqrt{\lambda}\,y) d\omega(\lambda) = \\
\quad \quad \quad = \int\limits_{-\infty}^{+\infty} \cos(\sqrt{\lambda}\,x) \cos(\sqrt{\lambda}\,y) d\omega(\lambda) +\int\limits_{-\infty}^{+\infty}  \int\limits_0^x K(x,t) \cos(\sqrt{\lambda}\,t) dt \cos(\sqrt{\lambda}\,y) d\omega(\lambda) = \\
\quad \quad \quad = f(x,y) + \int\limits_0^x K(x,t) \int\limits_{-\infty}^{+\infty} \cos(\sqrt{\lambda}\,t) dt \cos(\sqrt{\lambda}\,y) d\omega(\lambda) = \\
\quad \quad \quad = f(x,y) + \int\limits_0^x K(x,t) f(t,y)dt = \\
\quad \quad \quad  = - K(x,y),
\end{array} \]
by the Gelfand-Levitan equation \eqref{eq:GelfandLevitan}. \qed
\begin{cor} The following formula holds
\begin{equation} \label{eq:Xphi}
 \mathbb X(x) \phi(x,\lambda) =  \cos(\sqrt{\lambda}\,x) = B^*(x) \bbmatrix{1\\0}
\end{equation}
\end{cor}
\noindent\textbf{Proof:} Immediate from \eqref{eq:XDefInt} and the formulas \eqref{eq:KxtBasic}, \eqref{eq:phicos}. \qed

It turns out that we can find the inverse of $\mathbb X(x)$ as follows.
\begin{defn} The operator $\mathbb Y(x):\mathcal H\rightarrow\mathcal H$ is defined as follows
\begin{equation} \label{eq:YDefInt}
 (\mathbb Y(x) g)(\lambda) = g(\lambda) - \int\limits_0^x \phi(y,\lambda)\int\limits_{\mathbb R} \phi(y,\mu) g(\mu) d\omega(\mu) dy.
\end{equation}
\end{defn}
Notice that if we change the variables in $K(x,y)$ we will obtain a Kernel, which satisfies the defining initial value PDE for 
$K_1(x,y)$. By the uniqueness of the solution of this problem it follows that $K_1(x,y) = K(y,x)$. This enables us to show that $\mathbb Y(x)$ is
the inverse of the operator $\mathbb X(x)$.
Similarly to the calculation of $\mathbb X(x) \phi(x,\lambda)$ we obtain, using \eqref{eq:cosphi} that
\begin{equation} \label{eq:Ycos}
\mathbb Y(x) \cos(\sqrt{\lambda}\,x) = \mathbb Y(x) B^*(x) \bbmatrix{1\\0} = \phi(x,\lambda).
\end{equation}
Since as we mentioned earlier $[\phi(x,\lambda)]^*$ is a functional on $\mathcal H$ the following formula holds ($\mathbb Y^*(x)=\mathbb Y(x)$
from the definition)
\[ [\phi(x,\lambda)]^* = \bbmatrix{1&0} B(x)\mathbb Y(x).
\]
Moreover, using this notation, we also obtain that
\[ \dfrac{d}{dx} [\mathbb Y(x)] = - \phi(x,\lambda) [\phi(x,\lambda)]^*.
\]
Let us perform a simple calculation for the derivative of $\mathbb X(x)\mathbb Y(x)$:
\[ \begin{array}{lll}
\dfrac{d}{dx} [\mathbb X(x)\mathbb Y(x)] & = \dfrac{d}{dx} [\mathbb X(x)]\mathbb Y(x)] + \mathbb X(x)\dfrac{d}{dx} [\mathbb Y(x)] \\
& = B(x)\bbmatrix{1\\0}\bbmatrix{1&0}  B^*(x)  \mathbb Y(x) - \mathbb X(x) \phi(x,\lambda) [\phi(x,\lambda)]^* \\
& = B(x)\bbmatrix{1\\0} [\phi(x,\lambda)]^* - B(x) \bbmatrix{1\\0} [\phi(x,\lambda)]^* = \\
& = 0.
\end{array} \]
And since $\mathbb X(0)\mathbb Y(0) = I I = I$ is the identity operator, we obtain that $\mathbb X(x)$ is onto, and since it is of the form $I+C$ for a compact
operator $C$, it is also invertible by a Fredholm theorem, with the inverse $\mathbb Y(x)$. We could alternatively prove that $\dfrac{d}{dx} [\mathbb Y(x)\mathbb X(x)] =0$
to deduce the same result.
\begin{thm} The operator $\mathbb X(x)$ is invertible for all $x\geq 0$ and its inverse is the operator $\mathbb Y(x)$, defined in \eqref{eq:YDefInt}.
\end{thm}

Let us take
\begin{equation} \label{eq:KDef}
K'(x,y) = -\bbmatrix{1&0} B^*(x) \mathbb X^{-1}(x) B(y) \bbmatrix{1\\0}.
\end{equation}
Then Gel'fand-Levitan equation \eqref{eq:GelfandLevitan} holds
\[ \begin{array}{llllll}
K'(x,y) + f(x,y) + \int\limits_{0}^x K(x,t) f(t,y) dt = \\
= K'(x,y) + f(x,y) - \int\limits_{0}^x \bbmatrix{1&0} B^*(x) \mathbb X^{-1}(x) B(t) \bbmatrix{1\\0}  \bbmatrix{1&0} B^*(t) \mathbb X_0^{-1} B(y) \bbmatrix{1\\0}dt = \\
= K'(x,y) + f(x,y) - \bbmatrix{1&0} B^*(x) \mathbb X^{-1}(x)  \int\limits_{0}^x B(t) \sigma_2 B^*(t) dt \mathbb X_0^{-1} B(y) \bbmatrix{1\\0} = \\
= \text{ using vessel condition \eqref{eq:DB} } = \\
= K'(x,y) + f(x,y) - \bbmatrix{1&0} B^*(x) \mathbb X^{-1}(x) (\mathbb X(x) - \mathbb X_0) \mathbb X_0^{-1} B(y) \bbmatrix{1\\0} = \\
= K'(x,y) + f(x,y) - \bbmatrix{1&0} B^*(x) \mathbb X_0^{-1} B(y) \bbmatrix{1\\0} +
\bbmatrix{1&0} B^*(x) \mathbb X^{-1}(x)  B(y) \bbmatrix{1\\0} = 0.
\end{array} \]
By the uniqueness of the solution of \eqref{eq:GelfandLevitan} we conclude that $K'(x,y)=K(x,y)$ and the potential $q(x)$,  which is obtained from $K(x,x)$ using \eqref{eq:dKxxGL} is
\begin{equation} \label{eq:q=q'}
 q(x) = 2 \dfrac{dK(x,x)}{dx} = 2 \dfrac{dK'(x,x)}{dx} = -2 \dfrac{d}{dx} \bbmatrix{1&0} B^*(x) \mathbb X^{-1}(x) B(x) \bbmatrix{1\\0} =- 2  \dfrac{d^2}{dx^2} [\ln(\tau(x))],
\end{equation}
which coincides with \eqref{eq:qbeta}. In other words, we obtain the following
\begin{thm}[Realizing a potential on the half line] \label{thm:VesselFromq}
Let $\rho(\lambda)$ be the spectral function of the potential $q(x)$ for which $f(x,y)$ is defined using \eqref{eq:Deff'}.
Define a vessel
\[ \mathfrak{V} = (A, B(x), \mathbb X(x); \sigma_1, \sigma_2, \gamma, \gamma_*(x);
L^2_\omega,\mathbb C^2; [0,\infty)),
\]
obtained by the vessel construction, applied to $A = i\lambda$, $\mathbb X_0 = I$, $B_0 = \bbmatrix{1&0}$. Then this vessel realizes
the given potential $q(x)$.
\end{thm}
\noindent\textbf{Proof:} the differential equations \eqref{eq:DB}, \eqref{eq:DX} are satisfied by the construction. The Lyapunov equation \eqref{eq:Lyapunov}
for $x=0$ is satisfies, since then $A \mathbb X_0 + \mathbb X_0 A^* = 0$ for $\mathbb X_0=I$ and $A$ - anti selfadjoint by its definition. Notice also that
$B_0 = \bbmatrix{1\\0}$, so that $B_0 \sigma_1 B_0^* =0$ and as a result \eqref{eq:Lyapunov} holds. By Lemma \ref{lemma:Redund} it will hold for all $x$ by the
permanency of the Lyapunov equation. The operator $\mathbb X(x)$ is invertible for all $x\in[0,\infty)$ 
with the inverse $\mathbb Y(x)$, defined in \eqref{eq:YDefInt} by the arguments proceeding the
Theorem. The potential of the vessel coincides with the given potential $q(x)$ by \eqref{eq:q=q'}.
\qed

\section{Generalized scattering theory on a line}
A closer examination of the previous results, which can also be seen from the GL theory, implies that the operators $B(x), \mathbb X(x)$ are globally defined for all 
$x\in\mathbb R$. Indeed the formula \eqref{eq:DefBdomega}
\[  B(x) = \bbmatrix{\cos(\sqrt{\mu} x) & i \sqrt{\mu}\sin(\sqrt{\mu}x)}
\]
and \eqref{eq:XDefInt}
\[ \mathbb X(x) = I + \int\limits_0^x B(y)\sigma_2B^*(y) dy \]
define operators for all $x\in\mathbb R$. This implies the following definition
\begin{defn} \label{def:PreVessel} The collection of operators and spaces 
\begin{equation} \label{eq:DefPreV}
\mathfrak{preV} = (A, B(x), \mathbb X(x); \sigma_1, \sigma_2, \gamma;\mathcal{H},\mathbb C^p),
\end{equation}
is called a \textbf{pre-vessel}, if $\mathbb X(x), A:\mathcal H\rightarrow\mathcal H$, 
$B(x):\mathbb C^p\rightarrow\mathcal H$ are differentiable linear operators, subject to the following conditions:
\begin{enumerate}
	\item the operator $A$ has a dense domain $D(A)$ with a symmetric spectrum with respect to the imaginary axis.
		$A$ is the generator of a $C_0$ semi-group on $\mathcal H$,
	\item $B(x)$ satisfies regularity assumptions \eqref{eq:ResBsigma2}, \eqref{eq:ResBgamma}
	 and the equation \eqref{eq:DB} for all $x\in\mathbb R$,
	\item $\mathbb X(x)$ is bounded, self adjoint for all $x\in\mathbb R$ and satisfies \eqref{eq:DX} and there exists a point $x_0\in\mathbb R$, in which
	$\mathbb X(x_0)$ is invertible,
	\item the Lyapunov equation \eqref{eq:Lyapunov} holds for all $x\in\mathbb R$, $\lambda\not\in\SPEC(A)$.
\end{enumerate}
\end{defn}
Tau function $\tau(x) = \det(\mathbb X^{-1}(x_0)\mathbb X(x))$ is well defined now and exists on $\mathbb R$. Moreover, if $\tau(a')=\tau(b')=0$ and
$\tau(x)\neq 0$ for $a'<x<b'$, for arbitrary $\mathrm I=[a,b]$, where $a'<a<b<b'$ we obtain a vessel on $\mathrm I$. As a result, we can make the following definition
\begin{defn} The \textbf{generalized potential} of a pre-vessel $\mathfrak{preV}$ is the function $\gamma_*(x)$, defined by \eqref{eq:Linkage} for all points where
$\tau(x)\neq 0$.
\end{defn}
As we can see for each pre-vessel $\mathfrak{preV}$ there is a \textit{unique} generalized potential, realized by it. So, a very important question arises, is whether
each generalized potential $\gamma_*(x)$ arises in this form. It is also important to classify the classes of such generalized potentials and to relate them to
classes of pre-vessels.
\begin{problem}[Generalized scattering theory] Characterize the class $\GP$ of generalized potentials, arising from pre-vessels.
\end{problem}
\begin{problem}[Generalized inverse scattering] Given a generalized potential $q(x)\in\GP$, find a pre-vessel realizing it.
\end{problem}
For example, using Gelfand-Levitan theory we obtain:
\begin{thm} Let $q(x) (0\leq x<\infty)$ be a potential for which the function $f(x,y)$ possesses the form \eqref{eq:Deff'}. Then there
exists a pre-vessel, whose generalized potential coincides with $q(x)$ on the positive real line.
\end{thm}

\section{\label{sec:KdV}Solution of the KdV equation in the first quadrant}
Let as consider the initial value problem of \eqref{eq:KdV}.
Suppose that for a given $q(x)$ we have
constructed a vessel $\mathfrak{V}$ as in Theorem \ref{thm:VesselFromq}. Let us evolve the operators of the vessel $\mathfrak{V}$, using the following conditions:
\begin{eqnarray} 
\label{eq:DBt}
 \dfrac{\partial}{\partial t} B(x,t) = i A \dfrac{\partial}{\partial x} B(x,t), \quad B(x,0)=B(x) \\
\label{eq:DXt} \dfrac{\partial}{\partial t} \mathbb X = 
 i A B \sigma_2 B^* - i B\sigma_2 B^* A^* +
 i B\gamma B^*, \quad \mathbb X(x,0)=\mathbb X(x).
\end{eqnarray}
For these equations to hold we have to add some regularity assumptions:
\begin{equation} \label{eq:KdVReg}
 B(x,t)\sigma_2\mathbb C^p\subseteq D(A), \quad B(x,t)\gamma\mathbb C^p\subseteq\mathcal H. 
\end{equation}
Notice that since the operator $A=i\lambda$ is ``diagonal'' in our setting the equation \eqref{eq:DBt} is a sort of wave
equation, whose solution is very easy to create:
\begin{equation} \label{eq:BEvolved}
B(x,t) = \bbmatrix{\cos(\sqrt{\lambda} x - (\sqrt{\lambda})^3t) & i \sqrt{\lambda}\sin(\sqrt{\lambda}x- (\sqrt{\lambda})^3t)}
\end{equation}
Here $x,t\geq 0$ and $\cos(\sqrt{\lambda}) = \cosh(\sqrt{|\lambda|})$ if $\lambda<0$ as usual in 
Gelfand-Levitan notations. 
\begin{thm} Let $B(x,t)$ be a solution of the partial differential equations \eqref{eq:DB}, \eqref{eq:DBt}.
Then the formula for $\mathbb X(x,t)$ is as follows
\begin{multline} \label{eq:XEvolved}
\mathbb X(x,t) = \mathbb X_0 + \int_0^x B(y,t)\sigma_2B^*(x,t) dy + \\
+ \int_0^t [i A B(0,s) \sigma_2 B^*(0,s) - i B(0,s)\sigma_2 B^*(0,s) A^* +
 i B(0,s)\gamma B^*(0,s)] ds
\end{multline}
\end{thm}
\noindent\textbf{Proof:} The first equation \eqref{eq:DX} is immediate, since the last term depends on $t$ only.
In order to prove \eqref{eq:DXt}, we can calculate using \eqref{eq:DB} and \eqref{eq:DBt} that 
\[ \dfrac{\partial}{\partial t} [B(y,t)\sigma_2B^*(y,t)] = 
\dfrac{\partial}{\partial y} [i A B(y,t) \sigma_2 B^*(y,t) - i B(y,t)\sigma_2 B^*(y,t) A^* +
 i B(y,t)\gamma B^*(y,t)]
\]
and the result follows by applying $\int\limits_{0}^x dy$ to both sides. \qed
\begin{defn}\label{def:KdVVessel} The KdV evolutionary vessels is defined as follows
\[ \mathfrak{V} = (A, B(x,t), \mathbb X(x,t); \sigma_1, \sigma_2, \gamma, \gamma_*(x,t);
L^2_\omega,\mathbb C^2; \mathrm I),
\]
where $A=i\mu$ and $B(x,t), \mathbb X(x,t)$ are defined by \eqref{eq:BEvolved}, \eqref{eq:XEvolved} and satisfy the
regularity assumptions \eqref{eq:KdVReg}.
\end{defn}
Notice that for each fixed $t$, we obtain a vessel, for which all operators depend on $t$ and realize a potential $q(x,t)$. For the KdV equation \eqref{eq:KdV}
we have to require the existence of the third derivative with respect to $x$. In this case, it turns out 
that the regularity assumptions \eqref{eq:KdVReg} imply that the expressions
\[ \bbmatrix{1&0} B^*(x,t) \mathbb X^{-1}(x,t)A^i B(x,t)\bbmatrix{1\\0}, \quad i=0,1
\]
are meaningful. This enables to use the theory of evolutionary vessels, or more precisely
\cite[Theorem 10]{bib:KdVVessels}:
\begin{thm}[Vessel=Global Solution in First Quadrant]\label{thm:KdVequation}
Let $\mathfrak{V}_{KdV}$ be a KdV vessel, realizing a given potential $q(x)$, defined on the half line with 4 continuous derivatives. 
Then the potential $q(x,t)$ of the output SL equation \eqref{eq:OutCC} satisfies KdV equation \eqref{eq:KdV}:
\[ q'_t = - \dfrac{3}{2} q q'_x + \dfrac{1}{4} q'''_{xxx}
\]
for all $x,t\geq 0$.
\end{thm}

Let us calculate the value of $q(0,t), t\geq 0$. Using formula \eqref{eq:q=q'} we obtain that
\[ q(0,t) = -2 \dfrac{\partial}{\partial x} \bbmatrix{1&0}B^*(x,t)\mathbb X^{-1}(x,t) B(x) \bbmatrix{1\\0}|_{x=0}.
\]
Using \eqref{eq:DX-1BBounded} and \eqref{eq:DB} it becomes
\[ -\dfrac{1}{2} q(0,t) = \TR\bbmatrix{0&i\\-i&0}B^*(x,t)\mathbb X^{-1}(x,t) B(x))|_{x=0} - (\dfrac{\tau'}{\tau})^2|_{x=0}.
\]
Since $\dfrac{\tau'}{\tau} = \int_0^x q(y,t)dy$ at zero this term vanishes and we obtain
\[ -\dfrac{1}{2} q(0,t) = \TR\bbmatrix{0&i\\-i&0}B^*(x,t)\mathbb X^{-1}(x,t) B(x,t))|_{x=0} 
\]
Moreover, from \eqref{eq:Ycos} and \eqref{eq:DX-1BBounded} it follows that
\[ B^*(x,t)\mathbb X^{-1}_{x=0} = [\bbmatrix{\phi(0,t)\\-i(\phi'(0,t)+\dfrac{\tau'(0,t)}{\tau(0,t)}\phi(0,t)}]^* = [\bbmatrix{1\\0}]^*,
\]
and as a result
\[ \begin{array}{lll}
-\dfrac{1}{2} q(0,t) & = (\TR\bbmatrix{0&i\\-i&0}B^*(x,t)\mathbb X^{-1}(x,t) B(x))|_{x=0} =
\TR\bbmatrix{0&i\\-i&0}[\bbmatrix{1\\0}]^* B(0,t) = \\
& = \TR( [\bbmatrix{1\\0}]^*\bbmatrix{\cos((\sqrt{\lambda})^3t) & -i \sqrt{\lambda}\sin((\sqrt{\lambda})^3t)}\bbmatrix{0&i\\-i&0}\\
& = - \int_{\mathbb R} \sqrt{\lambda}\sin((\sqrt{\lambda})^3t) d\omega(\lambda).
\end{array} \]
in other words, it is assigned by the vessel construction:
\begin{equation} \label{eq:q0t}
q(0,t) = 2 \int_{\mathbb R} \sqrt{\lambda}\sin((\sqrt{\lambda})^3t) d\omega(\lambda).
\end{equation}
In order to satisfying different initial conditions for $q(0,t)$ we have to perturb the
original spectral function $\omega(\lambda)$ as follows. Suppose that the measure $\nu(\lambda)$ satisfies:
\[ \int\limits_{-\infty}^{+\infty}  \cos(\sqrt{\lambda}\,x) \cos(\sqrt{\lambda}\,y) d\nu(\lambda) = 0, \quad \forall x,y\geq 0.
\]
Let us define a new measure
\[ \rho'(\lambda) = \rho(\lambda) + \nu(\lambda), \quad \omega'=\omega+\nu.
\]
Then the formula \eqref{eq:Deff'} with $\rho'(\lambda)$ will produce the same $f(x,y)$ and hence by the uniqueness of GL equation, the same $K(x,y)$, and hence the same potential.
But the formula for the initial value for $q(0,t)$ will become
\[ q'(0,t) = 2 \int_{\mathbb R} \sqrt{\lambda}\sin((\sqrt{\lambda})^3t) d\omega'(\lambda) = 
q(0,t) + 2 \int_{\mathbb R} \sqrt{\lambda}\sin((\sqrt{\lambda})^3t) d\nu(\lambda),
\]
creating a perturbation of the original potential.

Finally we notice that the uniqueness of the solution of the KdV equation in this case ``follows'' from the fact that choosing a spectral function $\rho(\lambda)$
we actually uniquely determine a generalized vessel of the corresponding pre-vessel. As a result, the solution can be expected to be unique, as it is parallel
to the solution of the KdV equation on the whole line.

Finally, we obtain that evolving the operators as in the KdV case, we can create solutions of the KdV equation for generalized potentials using of pre-vessels.
\begin{defn} The KdV evolutionary pre-vessels is defined as follows
\[ \mathfrak{preV} = (A, B(x,t), \mathbb X(x,t); \sigma_1, \sigma_2, \gamma;
L^2_\omega,\mathbb C^2),
\]
where $A=i\mu$ and $B(x,t), \mathbb X(x,t)$ are defined by \eqref{eq:BEvolved}, \eqref{eq:XEvolved} and satisfy the
regularity assumption \eqref{eq:KdVReg}.
\end{defn}
\begin{thm} Given a KdV pre-vessel, the generalized potential $\gamma_*(x,t)$ exists for values $x,t\in\mathbb R^2$
for which $\tau(x,t)\neq 0$. At this point KdV equation \eqref{eq:KdV} holds.
\end{thm}

\section{\label{sec:supres}Supplementary results}
\subsection{Models of vessels in the case $var[\rho(\lambda)]<\infty$}
In the case $var[\rho(\lambda)]<\infty$ it is known that the spectral function can be represented as follows
\[ d\rho(\lambda) = f'_c(\lambda)d\lambda + \sum_{n=0}^\infty b_n \delta(\lambda-\lambda_n),
\]
where $f'_c(\lambda)$ is an integrable function on $\mathbb R$ and $b_n \delta(\lambda-\lambda_n)$ is the usual point measure, concentrated at $\lambda_n$
with the value $b_n$ there. We construct $\omega(\lambda)$ from this function
as at the Background
\[ \omega(\lambda) = \left\{ \begin{array}{ll}
\rho(\lambda) - \dfrac{2}{\pi}\sqrt{\lambda}, & \lambda > 0, \\
\rho(\lambda), & \lambda<0
\end{array}\right.
\]
so that $a(x) = \int_{\mathbb R} \dfrac{\cos(\sqrt{\lambda}x)}{\lambda} d\omega(\lambda)$ possesses fourth continuous derivative.

Define an element $\bbmatrix{g(\lambda)\\ g_n }\in\mathcal H$ to be a vector consisting of a function $g(\lambda)$ and an infinite sequence $\{g_n\}$ as follows:
\[ \mathcal H = \left\{ \bbmatrix{g(\lambda)\\ \{g_n\} } \mid \int\limits_{\mathbb R} |g(\lambda)|^2 f_c(\lambda) d\lambda + \sum_{n=1}^\infty |g_n|^2 b_n < \infty\right\}.
\]
Define next
\[ \begin{array}{lll}
A \bbmatrix{g(\lambda)\\ \{g_n\} } & = \bbmatrix{i\lambda g(\lambda)\\ \{ i\lambda_n g_n\} }, \\
\mathbb X_0 & =I, \\
B_0 & = \bbmatrix{ \bbmatrix{1 & 0} \\ \{ \bbmatrix{ 1 & 0} \} }
\end{array} \]
And apply the standard construction of a vessel to this data. Then
\[ B(x) = \bbmatrix{ \bbmatrix{\cos(\sqrt{\lambda}x) & i\sqrt{\lambda}\sin(\sqrt{\lambda}x)} \\
				\{ \bbmatrix{\cos(\sqrt{\lambda_n}x) & i\sqrt{\lambda_n}\sin(\sqrt{\lambda_n}x) \} } }
\]
and
\[ \mathbb X(x) = I + \int_0^x B(y)\sigma_2 B^*(y)dy
\]
will construct a vessel $\mathfrak V$ using Theorem \ref{thm:VesselFromq}. Evolving of this vessel with respect to $t$ will produce
\[ B(x,t) = \bbmatrix{ \bbmatrix{\cos(\sqrt{\lambda}x- (\sqrt{\lambda})^3t) & i\sqrt{\lambda}\sin(\sqrt{\lambda}x -(\sqrt{\lambda})^3t )} \\
				\{ \bbmatrix{\cos(\sqrt{\lambda_n}x(\sqrt{\lambda_n})^3t) & i\sqrt{\lambda_n}\sin(\sqrt{\lambda_n}x(\sqrt{\lambda_n})^3t) \} } }
\]
and the corresponding $\mathbb X(x,t)$, defined by \eqref{eq:XEvolved} producing a solution of the KdV equation \eqref{eq:KdV} with the initial $q(x)$
possessing the spectral function $\rho(\lambda)$ ($var[\rho(\lambda)]<\infty$).

Notice that it can be generalized to a model involving (tempered) distributions, in which case the formula \eqref{eq:Deff'} will hold at the sense of
distributions. This is a separate from this work project.

\subsection{Construction of a vessel using Fadeyev theory}
The same approach can implemented for the setting of Fadeyev theory \cite{bib:FaddeyevII}. For this it is enough to
consider Jost solutions $\phi(x,\lambda)$ and $\psi(x,\lambda)$ satisfying
\begin{eqnarray*}
\phi(s,x) = \cos(sx) + \int_0^x \dfrac{\sin(s(x-y))}{s} q(y) \phi(s,y) dy \\
\psi(s,x) = \dfrac{\sin(sx)}{s} + \int_0^x \dfrac{\sin(s(x-y))}{s} q(y) \psi(s,y) dy,
\end{eqnarray*}
for each $s$ at the upper half plane. this actually means that these two functions satisfy \eqref{eq:SL} with the given $q(x)$ and the initial conditions
\[ \phi(s,0)=1, \phi'(s,0)=0, \quad \psi(s,0)=0, \psi'(s,0)=1.
\]
Then define the following function, denoting $\beta(x)=\int_0^x q(y)dy$
\[ S(\lambda,x)=\bbmatrix{\phi(s,x)&i\psi(s,x)\\-i(\phi'(s,x)-\beta(x)\phi(s,x)) & \psi'(s,x)-\beta(x)\psi(s,x)} \bbmatrix{\cos(sx)&i\sin(sx)\\ s\sin(sx) &\cos(sx)}^{-1},
\]
where we notice that all the functions actually depend on $s^2=-i\lambda$.
Using the Fadeyev's work \cite{bib:FaddeyevII}, or more precisely the work of Agranovich-Marchenko \cite{bib:AgMarch} it follows that this function is well defined for all $x$ and satisfies the following differential equation
\[ S'(\lambda,x) = \sigma_1^{-1} (\sigma_1\lambda+\gamma_*(x)) S(\lambda,x) - S(\lambda,x) \sigma_1^{-1} (\sigma_1\lambda+\gamma),
\]
for $\gamma_*(x)$, defined by
\[ \gamma_*(x) = \bbmatrix{-i(\beta'(x)-\beta^2(x)) & -\beta(x)\\\beta(x)&i}.
\]
At $x=0$ it is a function of the complex variable $\lambda$, which has jumps on the positive part of the imaginary axis and a finite number of isolated points
on the negative imaginary part. In \cite{bib:GenVessel} it is shown that there exists a vessel with this transfer function with initial $S(\lambda,0)$,
where the operator $A$ is again diagonal. Applying KdV vessel construction for this $S(\lambda,0)$, we will obtain a solution of \eqref{eq:KdV},
analogous to the classical solution using Fadeyev's inverse scattering. Notice that in this case the measure $\omega(\lambda)$ satisfies $var[\omega(\lambda)]<\infty$
and we can use the standard construction to model the vessel.

\subsection{Examples of equations possessing vessel solutions}
The work in this section is presented in \cite{bib:ENLS} for the case of bounded operators, but the ideas presented here
are immediatly generalizable for an unbounded $A$, generating a $C_0$ semi-group.
\begin{defn} NLS vessel parameters are defined as follows:
\[ \sigma_1 = \bbmatrix{1&0\\0&1},\quad
\sigma_2 = \dfrac{1}{2} \bbmatrix{1&0\\0&-1},
\quad \gamma =\bbmatrix{0 & 0 \\0 & 0},
\]
\end{defn}
And constructing a vessel in a similar to SL manner from a realized function $S(\lambda,0)$, we will obtain solutions of the NLS equation:
\[
 \dfrac{\partial}{\partial x} y(\lambda,x) - (- i \lambda) \sigma_2 (\lambda,x) + \bbmatrix{0&\beta(x)\\-\beta^*(x)&0} u(x,\lambda) = 0,
\]
where $\beta(x)$ is easily derived from the vessel.
Evolving next the obtained vessel using \eqref{eq:DBt}, \eqref{eq:DXt}, we will obtain solutions of the evolutionary NLS:
\begin{equation} \label{eq:ENLS}
i \beta_t + \beta_{xx} + 2 |\beta|^2 \beta = 0, \quad \beta(x,0) = \beta(x),
\end{equation}

For the Canonical systems, one can use the following parameters:
\begin{defn} Canonical systems vessel parameters are defined as follows:
\[ \sigma_1 = \bbmatrix{0&i\\-i&0},\quad
\sigma_2 = \bbmatrix{1&0\\0&1},
\quad \gamma =\bbmatrix{0 & 0 \\0 & 0},
\]
\end{defn}
and to use a similar theory of evolutionary vessels to create solutions of
\begin{equation} \label{eq:CanSysE}
\dfrac{\partial}{\partial t}\big( \sqrt{\dfrac{1}{(x+K)^2} - 4\beta^2}\big) = -2 \dfrac{\beta}{(x+K)^2} + \beta_{xx},
\end{equation}
where it turns out that for afixed $K\in\mathbb R$ 
$\gamma_*(x,t)=\bbmatrix{-2 i \beta(x,t)&i h(x,t)\\ih(x,t)&2i\beta(x,t)}$ with $h(x,t)=\sqrt{\dfrac{1}{(x+K)^2} - 4\beta^2}$.

\bibliographystyle{alpha}
\bibliography{../../biblio}

\end{document}